\documentclass[11pt, a4paper]{article}
\usepackage{amsmath}
\usepackage{amsfonts}
\usepackage{bm}
\newcommand\BBE{{\mathbb {E}}}

\newcommand\R{{\mathbb {R}}}
\newcommand\N{{\mathbb {N}}}

\newtheorem{Theorem}{Theorem}
\newtheorem{corollary}[Theorem]{Corollary}

\newtheorem{Lemma}[Theorem]{Lemma}

\newtheorem{Proposition}[Theorem]{Proposition}
\newtheorem{remark}[Theorem]{Remark}

\renewcommand{\P}{\ensuremath{\mathbb {P}}}
\newcommand{\E}{\mathbb{E}}

\newcommand{\p}{ { \mathbb P} }

\newcommand\beq{\begin{equation}}
\newcommand\eeq{\end{equation}}

\begin{document}

\title{Strong approximations for a class of dependent random variables with semi exponential tails}

\author{C. Cuny\footnote{Christophe Cuny, Univ Brest, UMR 6205 CNRS 6205, LMBA, 6 avenue Victor Le Gorgeu, 29238 Brest}, J. Dedecker\footnote{J\'er\^ome Dedecker, Universit\'e Paris Cit\'e, CNRS, MAP5, UMR 8145,
45 rue des  Saints-P\`eres,
F-75006 Paris, France.}
and 
 F. Merlev\`ede \footnote{Florence Merlev\`ede, LAMA,  Univ Gustave Eiffel, Univ Paris Est Cr\'eteil, UMR 8050 CNRS,  \  F-77454 Marne-La-Vall\'ee, France.}}

\maketitle

\begin{abstract} 
We give rates of convergence in the almost sure invariance principle for sums of dependent random variables with semi exponential tails,  whose coupling coefficients decrease at a sub-exponential rate. We show that the rates in the strong invariance principle are in powers of $\log n$. We apply our results to iid products of random matrices.
\end{abstract}



\section{Introduction}

Let $(\Omega, {\mathcal A}, {\mathbb P})$ be a probability space, and let $(\varepsilon_i)_{i \geq 1}$ be independent and identically distributed (iid) random variables defined on $\Omega$, with values in a measurable space $G$ and with common distribution $\mu$. Let $W_0$ be a random variable defined on $\Omega$ with values in a measurable space $X$, independent of 
$(\varepsilon_i)_{i \geq 1}$, and let 
$F$ be a measurable function from $G\times X $ to $X$. For any $n \geq 1$, define
\[
W_{n} = F ( \varepsilon_{n}, W_{n-1} ) \, ,
\]
and assume that $(W_n, n \geq 1)$ has a stationary distribution $\nu$. 
Let now $h$ be a measurable function from $G\times X $ to ${\mathbb R}$ and define, for any $n \geq 1$,
\beq \label{functionofarandomiterates}
X_n = h ( \varepsilon_{n}, W_{n-1} ) \, .
\eeq
Then $(X_n)_{n \geq 1}$ forms a stationary sequence with stationary distribution, say  $\pi$.  Let $( {\mathcal G}_{i} )_{i \in {\mathbb Z}}$ be the non-decreasing filtration defined as follows: for any $i < 0$, 
${\mathcal G}_{i} =\{ \emptyset, 
\Omega\}$, ${\mathcal G}_{0} = \sigma (W_0)$ and for any $i \geq 1$, $ 
{\mathcal G}_{i} = \sigma ( \varepsilon_{i }, \ldots,  \varepsilon_{ 1 }, W_0 ) $. It follows that for any $n \geq 1$, $X_n$ is ${\mathcal G}_{n}$-measurable.  

We can also consider the following model
\beq \label{functionofiid}
X_n = h ( \varepsilon_{n}, \varepsilon_{n-1} , \ldots) \, ,
\eeq
which is in fact included in the preceding situation, by taking 
$W_0=( \varepsilon_{0}, \varepsilon_{-1} , \ldots)$. 

In what follows we assume that $h$ is such that $\E(X_n)=0$.  
For any $n \geq 0$, let  us define the sequence $(\delta (n) )_{n \geq 0}$ by
\[
\delta (0)=  \E (|X_1|) \,  \text{ and } \, \delta (n) = \Vert X_n - X_n^* \Vert_1 \, , \, n \geq 1 \, , 
\]
where $X_n^* $ is defined as follows :  
 let $W_0$ and $W_0^*$ be random variables with law $\nu$, and such that $W_0^*$ is independent of $(W_0, (\varepsilon_i)_{i \geq 1})$. For any $n \geq 1$,  let 
\beq \label{defXnstartWnstar}
X^*_n = h ( \varepsilon_{n}, W^*_{n-1} )\,  \mbox{ with } \,  W^*_{n} = F ( \varepsilon_{n}, W^*_{n-1} ) \, .
\eeq
Note that the coefficients $\delta(n)$ are well defined   if $\pi$ has a moment of order $1$.

In the paper \cite{CDM}, we assumed that  the  coupling coefficients $\delta_n$ decrease at a polynomial rate,  and we used a variant of the Berkes-Liu-Wu method (see \cite{BLW14}) to establish strong approximation results for partial sums. Another more restrictive coefficient is also considered in \cite{CDM} (see also \cite{CDMbis} for variants of these coefficients and the resulting conditions). In \cite{CDM} we  applied our results to various classes of (possibly non irreducible) Markov chains; in particular, we obtained optimal rates in the almost sure invariance principle (ASIP) for some functions of the left random walk on $GL_d({\mathbb R})$ under polynomial moment conditions.  In \cite{CDKM0} and \cite{CDKM} we  adapted the proofs of these results to Birkhoff sums of H\"older observables  of non-uniformly expanding dynamical systems, using a representation  as functions of the trajectory of a particular Markov Chain due to Korepanov \cite{Kor}.  In the paper \cite{CDKM}, we consider  dynamical systems with exponential or sub-exponential decay of correlations, and we show that the rates in the ASIP are in powers of $\log n$. 

In the present paper, we consider the context of \cite{CDM}, assuming that the coefficients $\delta_n$ decrease at an exponential, sub-exponential or super-exponential rate (described by an index $\gamma_1$, see \eqref{conddelta}), and that the variables have an exponential, sub-exponential or super-exponential moments (described by an index  $\gamma_2$, see \eqref{condh}). Using arguments from \cite{CDKM}, we show that the rates  in the ASIP are in powers of $\log n$ (the exponent of the logarithm depending on  $\gamma_1$ and $\gamma_2$). We apply our results to the left random walk on $GL_d({\mathbb R})$, to iid products of positive matrices, and to a class of non-uniformly contracting auto-regressive processes.

\section{Main results}
We  assume that there exist two positive constants $\gamma_1 $ and $c$ 
such that 
\beq \label{conddelta}
\delta (n) \leq  \exp ( - c n^{\gamma_1} )  \ \text{ for any positive integer $n$}\, , 
\eeq
and there are constant $b \in ]0, \infty [$ and $\gamma_2 \in ]0, \infty ]$ such that 
\beq \label{condh}
\p (|X_1|  >t )  \leq  \exp (1 - (t/b)^{\gamma_2} )\ \text{ for any positive  $t$} \, ,
 \eeq
Note that when $\gamma_2 = \infty$ \eqref{condh} means that $\Vert X_n \Vert_{\infty} \leq b$ a.s. for any integer $n$.   When $X_1$ satisfies \eqref{condh} we say either that it has a semi exponential tail of order $\gamma_2$ or that it admits a sub-exponential or super-exponential moment of order $\gamma_2$.

\begin{Theorem} \label{thmsubexpo} Let $(X_n)_{n \geq 1}$ be defined  by \eqref{functionofarandomiterates} or  \eqref{functionofiid}.  Let  $S_n = \sum_{k=1}^n X_k$. Assume  that \eqref{conddelta} and \eqref{condh} are satisfied. Then $n^{-1} \E \big ( S^2_n  \big ) \rightarrow \sigma^2$ as $n \rightarrow \infty$ and one can redefine $(X_n)_{n \geq 1}$ without changing its distribution on a (richer) probability space on which 
there exist iid random variables $(N_i)_{i \geq 1}$ with common distribution ${\mathcal N} (0, \sigma^2)$, such that, as $n \rightarrow \infty$, 
\beq \label{resultSA}
 S_n - \sum_{i=1}^n N_i = O( (\log n )^{1 + \lambda}) \, \text{ ${\mathbb P}$-a.s. where  $\lambda= \frac{1}{\gamma_1} + \frac{1}{\gamma_2} $.}
\eeq
\end{Theorem}

\noindent {\it Proof of Theorem \ref{thmsubexpo}.} Let $\alpha = 1+\lambda$. The proof uses similar arguments as those developed in the proof of Theorem 1.6 in \cite{CDKM} with the main difference that the random variables $X_k$ are no longer  bounded.  we then have to truncate them.  With this aim, let us first define for any positive integer $k$, 
\beq \label{defMk}
M_k= c_1 k^{1/\gamma_2}  \mbox{ where } c_1 =b (2 \log 3)^{1 / \gamma_2} \, , 
\eeq
(recall that $b$ is the constant appearing in condition \eqref{condh}) and
\[
\varphi_k (x) = \big ( x \wedge M_k \big ) \vee \big ( - M_k\big ) \, .
\]
Then we set 
\[
{ X}_{k,j} =\varphi_k (X_j)  - \E \varphi_k (X_j) \,  \mbox{ and }  \, 
W_{k, \ell} = \sum_{i=1+3^{k-1}}^{\ell +  3^{k-1}} { X}_{k,i}  \, .
\]
For $n \geq 2$, let $h_n =  \lceil (\log n)/(\log 3)\rceil $ (so that $h_n$ is the unique integer such that $3^{h_n -1} < n \leq 3^{h_n}$). Note that 
\[
S_n= \sum_{k=1}^{h_n-1} \sum_{i=1+3^{k-1}}^{ 3^k} X_i + \sum_{i=1+3^{h_n -1} }^n { X}_{i} \, , 
\]
and set 
\[
S_n^{\dag}= \sum_{k=1}^{h_n-1}   W_{k, 3^k - 3^{k-1}}+ \sum_{i=1+3^{h_n -1} }^n { X}_{h_n,i} \, . 
\]
Let also $S_1^{\dag}=0$.  
We first  prove that 
\beq \label{p1KMT}
\max_{1 \leq i \leq n} \big | S_i - S_i^{\dag} \big | =  O( (\log n)^{\alpha}) \, \text{ ${\mathbb P}$-a.s.}
\eeq
With this aim, it suffices to prove via the Kronecker lemma and stationarity that 
\beq \label{p1KMTp1}
 \sum_{k \geq 1} k^{- \alpha} 3^k \E ( | X_1 | -M_k)_+ < \infty \, .
\eeq
Obviously this holds when  $\gamma_2 = \infty$. Assume now that $ \gamma_2 < \infty$ and note that since $3^k = \exp \big ( \frac{1}{2} \Big (\frac{M_k}{b} \big )^{\gamma_2}\big )$, 
\[
3^k \E ( | X_1 | -M_k)_+ =3^k \int_{M_k}^{\infty} \p (|X_1|  >t ) dt  \leq \int_{M_k}^{\infty}  \exp (1 - 2^{-1}(t/b)^{\gamma_2} ) dt  \leq \frac{b 2^{1/\gamma_2}}{\gamma_2} \Gamma ( (1-\gamma_2)/\gamma_2 ) < \infty \, .
\]
Since $\alpha >1$, \eqref{p1KMTp1} is proved and so \eqref{p1KMT} also. 

Now, for any $k \geq 1$, let 
\beq \label{defmk}
m_k = [  c_2 k^{1/\gamma_1}]  +1 \, \text{ where $c_2 = \big ( 2 c^{-1} (\log 3) \big )^{1/\gamma_1}$} \, ,
\eeq
(recall that $c$ is the constant appearing in condition \eqref{conddelta}). We then define
\[
{\tilde X}_{k,j} = \E \big ( \varphi_k (X_j)   | \varepsilon_j, \varepsilon_{j-1}, \ldots, \varepsilon_{j-m_k}\big ) - \E \big ( \varphi_k (X_j) \big )  \text{ for any  $j \geq m_k$} \, \text{ and } \, {\widetilde W}_{k, \ell} = \sum_{i=1+3^{k-1}}^{\ell +  3^{k-1}} {\tilde X}_{k,i} \, .
\]
Finally, set  $ {\tilde S}_1 =0$ and for $n >1$, \[
 {\tilde S}_n= \sum_{k=1}^{h_n-1} {\widetilde W}_{k, 3^k - 3^{k-1}} +{\widetilde W}_{h_n, n- 3^{h_n -1}} \, . 
\]
Note that by Lemma 24 in \cite{CDM}, for any $j \geq m_k$, 
\beq \label{conslmaenlevertildedeltan}
\Vert X_{k,j}  - {\tilde X}_{k,j}  \Vert_1 \leq  \iint  \E (|  X_{m_k+1,x} - X_{m_k+1,y} |  ) \nu(dx) \nu (dy)  \leq \delta(m_k) . 
\eeq
Hence, for $k \geq k_0$ where $k_0$ is such that for any $\ell \geq k_0$,  $m_\ell \leq 3^{\ell-1}$, the upper bound \eqref{conslmaenlevertildedeltan} implies that
\[
\big \Vert  \max_{1 \leq \ell \leq 3^k - 3^{k-1}} \big |   { W}_{k, \ell} -  {\widetilde W}_{k, \ell} \big | \big \Vert_1 \leq  \sum_{i=1+3^{k-1}}^{ 3^{k}} \Vert X_{k,i}  - {\tilde X}_{k,i}  \Vert_1 \leq   3^k \delta (m_k) \ll 3^k \exp (-  c \times c_2^{\gamma_1}k)   \, .
\]
Hence, since $c \times c_2^{\gamma_1}  \geq  \log 3 $, 
\[
\sum_{k \geq k_0}  k^{-\alpha}  \big \Vert  \max_{1 \leq \ell \leq 3^k - 3^{k-1}} \big |   { W}_{k, \ell} -  {\widetilde W}_{k, \ell} \big | \big \Vert_1 < \infty \, ,
\]
which entails by  the Kronecker lemma that 
\beq \label{p2KMT}
\max_{1 \leq i \leq n} \big |  S_i^{\dag}  - {\tilde S}_i \big | =  o( (\log n )^{\alpha}) \, \text{ ${\mathbb P}$-a.s.}
\eeq
From \eqref{p1KMT}  and  \eqref{p2KMT},  the theorem will follow if one can prove that there exists a standard Brownian motion $B$ such that 
such that, as $n \rightarrow \infty$, 
\beq \label{resultinter}
\max_{1 \leq i \leq n} \big | {\tilde S}_i - B ( i \sigma^2)   \big | = O( (\log n )^{\alpha}) \, \text{ ${\mathbb P}$-a.s.}
\eeq
To prove \eqref{resultinter}, we do exactly as in Steps 3 and 4 of the proof of \cite[Theorem 1.6]{CDKM} by replacing the upper bound $\Vert \psi \Vert_{\infty}$ by 
$2 M_k$ (which a bound of the sup norm of the  truncated random variables ${ X}_{k,j} $), and by noticing that $M_k m_k \ll k^{1/\gamma}$. This leads to  (3.20) in \cite{CDKM} which combined with  (3.9) in \cite{CDKM} gives the following strong approximation: 
there exists a standard Brownian motion $B$ such that 
such that, as $n \rightarrow \infty$, 
\beq \label{resultinter2}
\max_{1 \leq i \leq n} \big | {\tilde S}_i - B (  \sigma_i^2)   \big | = O( (\log n )^{\alpha}) \, \text{ ${\mathbb P}$-a.s.}
\eeq
where $\sigma^2_i$ is defined by (3.21) in \cite{CDKM}.   It remains to identify the variance of the Brownian motion and to show that 
\beq \label{resultinter3}
\max_{1 \leq i \leq n} \big | B (  \sigma_i^2) - B( i \sigma^2)  \big | = O( (\log n )^{\alpha}) \, \text{ ${\mathbb P}$-a.s.}
\eeq
With this aim, we shall proceed as in step 3.4 in \cite{BLW14}  (see also Step 5 in \cite{CDKM}). We set 
 \beq \label{nuk1}
\nu_k =  \sigma^2 + 2 \sum_{i=-m_k}^{m_k} (   {\tilde \gamma}_{k,i} - \gamma_i ) - 2  \sum_{i \geq m_k+1} \gamma_i \, ,
\eeq
where \[
{\tilde \gamma}_{k,i}= {\rm cov} ({\tilde X}_{k,m_k} , {\tilde X}_{k,i+m_k} ) \, , \,  {\hat \gamma}_{k,i}= {\rm cov} ({ X}_{k,0} , { X}_{k,i} ) \, \text{ and } { \gamma}_{i} = {\rm cov} ({ X}_{0} , { X}_{i} ) \, .
\]
To prove \eqref{resultinter3}, we need to prove that 
\beq \label{cond1v_k}
( \log n)  \max_{ k \leq h_n} (m_k \nu_k)^{1/2} = O( (\log n )^{\alpha}  ) \, ,  
\eeq
and
\beq \label{cond2v_k}
3^k (  \nu_k^{1/2}  - \sigma)^2 = O( k^{2 \alpha} (\log k)^{-1}) \, .  
\eeq
By using the upper bound \eqref{conslmaenlevertildedeltan}, we have, for any $i \geq 0$, 
\begin{multline} \label{nuk2}
\big | {\tilde \gamma}_{k,i} - {\hat \gamma}_{k,i} \big |  =   \big |   {\rm cov} ({\tilde X}_{k,m_k} -  { X}_{k,m_k}, {\tilde X}_{k,i+m_k} )  +  {\rm cov} ( { X}_{k, m_k}, {\tilde X}_{k,i+m_k} -{X}_{k, i+m_k} )  \big |    \\ 
\leq   2 M_k  \Vert {\tilde X}_{k,m_k} -{ X}_{k,m_k} \Vert_1 +   2 M_k  \Vert  {\tilde X}_{k,i+m_k} - X_{k,i+m_k}   \Vert_1 \ll  k^{1/\gamma_2} \exp ( - c \times c_2^{\gamma_1} k )   \, .
\end{multline}
Next,  by using inequality (1.11a) in \cite{Ri00},  note that
\begin{multline*} 
\big | {\hat \gamma}_{k,i} - { \gamma}_{i} \big |  =   \big |   {\rm cov} ( X_{k,0} -  { X}_{0}, { X}_{k,i} )  +  {\rm cov} ( { X}_{0}, { X}_{k,i} -{X}_{i} )  \big |    \\ 
  \leq  2 \int_0^{1/2} Q_{|\varphi_k (X)|} (u) Q_{|g_k (X)|} (u) du  + 2 \int_0^{1/2} Q_{|X_1|} (u) Q_{|g_k (X)|} (u) du \leq  4   \int_0^1  Q (u) ( Q(u) - M_k )_+ du  \, ,
\end{multline*}
where $g_k (x) = (|x| -M_k)_{+}= x- \varphi_k (x)$.  But, for any positive  $t$, $
\p (|X_1|  >t )  \leq  \exp (1 - (t/b)^{\gamma_2} ): = G(t)$. Therefore, for any $u \in ]0, 1[$,
\[
Q(u) \leq G^{-1} (u) := b \big ( 1 - \log u \big )^{1/\gamma_2} \, .
\]
But
\[
G^{-1} (u)>M_k \iff  u^{-1} > e^{2 k \log 3 -1}  \, .
\]
So, overall, 
\begin{multline*} 
\big | {\hat \gamma}_{k,i} - { \gamma}_{i} \big |  \leq  4b^2   \int_0^1   \big ( \log ( e /u ) \big )^{2/\gamma_2}  {\bf 1}_{\{  u^{-1} > e^{2 k \log 3 -1}  \}  }du = 4 e b^2   \int^{\infty}_{e^{2 k \log 3 } }   \big ( \log  x  \big )^{2/\gamma_2}  x^{-2} dx \\
= 4 e b^2   \int^{\infty}_{{2 k \log 3 } }   x^{2/\gamma_2}  e^{-x} dx  \leq \kappa_{\gamma_2} b^2 \big ( 1 + k^{2/\gamma_2}  \big )  e^{-2 k \log 3 }  \, ,
\end{multline*}
where $\kappa_{\gamma_2} $ is a positive constant depending only on $\gamma_2$. 

Next, by using Proposition 1 in  \cite{DD03}, we derive
\[
\sum_{i\geq m_k} | {\rm cov} (X_0, X_i) | \leq  2  \sum_{i\geq m_k}  \int_0^{\delta(i)} Q \circ H^{-1} (u) du \leq  2  \sum_{i\geq m_k}  \int_0^{H ( \delta(i)) } Q^2 (u)  du   \, ,
\]
where  $H$ is the function defined on $[0,1]$ by $H( v) := \int_0^v Q(u) du $ and $H^{-1}$ is its inverse.  Since $Q(u) \leq  b \big ( 1 - \log u \big )^{1/\gamma_2}$, we infer that there exists a positive constant $\kappa $ depending on $c$, $b$, $\gamma_1$ and $\gamma_2$, such that, for any $i \geq 1$, 
\[
\int_0^{H ( \delta(i)) } Q^2 (u)  du \leq \kappa  i^{2 \gamma_1/\gamma_2} \exp (- c i^{\gamma_1} ) \, .
\]
Therefore, since $m_k =[c_2 k^{1/ \gamma_1} ] +1$, 
\[
\sum_{i\geq m_k} | \gamma_i | \leq {\tilde \kappa}    k^{ \frac{2}{ \gamma_2} + \frac{1}{ \gamma _1}  -1 }  \exp (- c \times c_2^{\gamma_1}  k ) \, ,
\]
where ${\tilde \kappa} $ is a positive constant depending only on $c$, $b$, $\gamma_1$,  $\gamma_2$ and $c_2$. So, overall, for any positive integer $k$, 
 \beq \label{nuk1maj}
\big | \nu_k - \sigma^2 \big |  \leq K  \times   k^{ \frac{2}{ \gamma_2} + \frac{1}{ \gamma _1} }  \Big \{  e^{-2 k \log 3 } +  \exp (- c \times c_2^{\gamma_1}  k )  \Big \}  \, .
\eeq
where $K$ is a positive constant depending only on $c$, $b$, $\gamma_1$,  $\gamma_2$ and $c_2$. 
This shows that $\nu_k \rightarrow \sigma^2 $, as $k \rightarrow \infty$. Hence \eqref{cond1v_k}  is satisfied since $( \log n)  \max_{ k \leq h_n} (m_k )^{1/2} \ll (\log n )^{1+ 1/(2 \gamma_1) } = O( (\log n )^{\alpha}  ) $ (indeed $\alpha = 1 + 1/\gamma_1  + 1/ \gamma _2$).  Note now that since $|\nu_k^{1/2} - \sigma |^2 \leq |\nu_k - \sigma^2|$,   \eqref{nuk1maj} implies \eqref{cond2v_k} since we have selected  $c_2 $ such that $ c \times c_2^{\gamma_1}  = 2 \log 3$. The proof is complete.  $\diamond$

\section{Applications}

\subsection{Products of iid invertible matrices \label{LRW}}
Let $(\varepsilon_n)_{n \geq 1}$ be independent random matrices taking values in $G= GL_d(\mathbb R)$, $d \geq 2$, with common distribution $\mu$.  
Let $A_0={\rm Id}$ and for every $n \geq 1$, $A_n = \varepsilon_n \cdots \varepsilon_1$.

Let $\Vert \cdot \Vert$ be the euclidean norm on ${\mathbb R}^d$ and for any $g \in G$, let $N(g) := \max ( \Vert g \Vert , \Vert g^{-1} \Vert)$ where  
$\|g\|=\sup_{\|x\|=1}
\|gx\|$. 
Recall that  $\mu $ has a moment of order $p \geq 1$ if  $\int_G (\log N(g) )^p \mu(dg) < \infty$ and a sub-exponential moment of order $\gamma \in ]0,1]$ if 
there exists $c>0$ such that 
\beq \label{condsubexpo}
\int_G {\rm e}^{c (\log N(g) )^{\gamma}}  \mu (d g) < \infty \, .
\eeq
Recall also that if $\mu$ admits a moment of order $1$ then 
\begin{equation}\label{LLN}
\lim_{n \rightarrow \infty} \frac{1}{n} \log \Vert A_n \Vert = \lambda_{\mu} \, \text{ ${\mathbb P}$-a.s.}, 
\end{equation}
where $\lambda_\mu:= \lim_{n\to +\infty} n^{-1}
\E(\log \|\varepsilon_n \cdots \varepsilon_1\|)$ is the so-called first Lyapunov exponent (see for instance \cite{FK}).

Let $X:=P_{d-1} ({\mathbb R}^d)$ be the projective space of ${\mathbb R}^d -\{0\}$ and write ${\bar x}$ as the projection of $x \in {\mathbb R}^d -\{0\}$ to $X$. $G$ is acting on $X$ as follows: $g\cdot \bar x=\overline{gx}$ for every $(g,x)\in G\times \R^d-\{0\}$.

\medskip

We assume that $\mu$ is strongly irreducible (i.e. that no proper finite union of subspaces of $\R^d$ are invariant by $\Gamma_\mu$, the closed semi-group generated by the support of $\mu$) and proximal (i.e. that there exists a matrix in $\Gamma_\mu$ 
admitting a unique (with multiplicity one) eigenvalue with maximum modulus). 
Under those assumptions (see e.g. Bougerol-Lacroix \cite{BL} or 
Benoist-Quint \cite{BQ}) it is well-known that there exists a unique invariant measure $\nu$ on ${\mathcal B} (X)$, meaning that for any continuous and bounded function $f$ from $X$ to $\mathbb R$,
\[
\int_X f(x) \nu(dx) = \int_G \int_X f( g \cdot x ) \mu(dg) \nu(dx) \, .
\]
The left random walk of law $\mu$ on $X$  is the process defined as follows. Let $W_0$ be an
$X$-valued  random variable 
 independent of $(\varepsilon_n)_{n \geq 1}$, and let $W_n =\varepsilon_n\cdot W_{n-1}$ for $n \geq 1$. Note that, if $W_0$ has distribution $\nu$, then 
$(W_n)_{n\geq 0}$ is a strictly stationary Markov chain. 

Our aim is to study the partial sums associated with the 
random sequence   $(X_n)_{n \geq 1}$ given by 
\[
X_n := h (\varepsilon_n, W_{n-1} ) \, , \, n \geq 1 \, , 
\]
where for every $g \in G$ and every ${\bar x} \in X$,
\[
 h ( g , {\bar x} ) = \log \Big ( \frac{\Vert g  x \Vert }{ \Vert x \Vert }\Big ) \, .
\]

As usual, we shall denote by $X_{n,{\bar x}}$ the random variable for which $W_0 = {\bar x}$. We then define $S_{n,{\bar x}}= \sum_{k=1}^n X_{n,{\bar x}}$. Note that,  for any $x \in S^{d-1}$, 
\begin{equation*}
S_{n, \bar x} = \sum_{k=1}^n X_{k, \bar x} = \log \Vert  A_n x \Vert \, .
\end{equation*}

Then, when we keep the notation $(X_n)_{n\ge 1}$, we have in mind that we are in stationary regime, i.e. that $W_0$ has law $\nu$. In this case we also define 
$$
S_n=\sum_{k=1}^n X_n =\log \|A_n V_0\|\, ,
$$
where $V_0$ is such that $\|V_0\|=1$ and $W_0=\bar V_0$.

\medskip

We denote by $\rho(g)$ the spectral radius of a matrix $g$. Applying Theorem  \ref{thmsubexpo}, the following strong  approximation with rate holds.

\begin{corollary} \label{KMTforcocycles} Let $\mu$ be a proximal and strongly irreducible probability measure on ${\mathcal B}(G)$. Assume that $\mu$ has a sub-exponential moment of order $\gamma \in ]0,1]$. Then $n^{-1} \E_{\nu} \big ( (S_n - n \lambda_{\mu} )^2\big ) \rightarrow \sigma^2$ as $n \rightarrow \infty$ and for every (fixed) ${\bar x} \in X$, one can redefine $(S_{n,{\bar x}})_{n \geq 1}$ without changing its distribution on a (richer) probability space on which 
there exist iid random variables $(N_i)_{i \geq 1}$ with common distribution ${\mathcal N} (0, \sigma^2)$, such that,
\begin{equation*} \label{asip-RW}
S_{n,{\bar x}} - n \lambda_{\mu} - \sum_{i=1}^n N_i = O( (\log n)^{1+2/\gamma}) \, \text{a.s.}
\end{equation*}
Moreover, the result remains true  if we replace the sequence $(S_{n,\bar x})_{n\ge 1}$ 
with $(S_n)_{n\ge 1}$, $(\log \|A_n\|)_{n\ge 1}$, 
$(\log \rho(A_n))_{n\ge 1}$ or $(\log |\langle A_n x, 
y\rangle |)_{n\ge 1}$,  for some $x,y\in S^{d-1}$.
\end{corollary}
\begin{remark} It follows from item $c)$ of Theorem 4.11 of 
Benoist-Quint \cite{BQ} that $\sigma>0$ if $\mu$ is strongly irreducible and the image of $\Gamma_\mu$ 
in $PGL_d(\R)$ is unbounded.  
\end{remark}
\begin{remark} \label{remsuperexpo}
If $\mu$ has a super-exponential moment meaning that \eqref{condsubexpo} holds with $\gamma >1$, we infer from the proof of Corollary \ref{KMTforcocycles} that for every (fixed) ${\bar x} \in X$, one can redefine $(S_{n,{\bar x}})_{n \geq 1}$ without changing its distribution on a (richer) probability space on which 
there exist iid random variables $(N_i)_{i \geq 1}$ with common distribution ${\mathcal N} (0, \sigma^2)$, such that,
\[
 S_{n,{\bar x}} - n \lambda_{\mu} - \sum_{i=1}^n N_i = O( (\log n)^{2+1/\gamma}) \, \text{a.s.}
\]
In particular if $\mu$ has compact support, the rate in the almost sure invariance principle is of order  $O( (\log n)^{2})$. Again, this extension holds with $(S_n)_{n\ge 1}$, $(\log \|A_n\|)_{n\ge 1}$, 
$(\log \rho(A_n))_{n\ge 1}$ or $(\log |\langle A_n x, 
y\rangle|)_{n\ge 1}$,  for some $x,y\in S^{d-1}$.
\end{remark}

\noindent \emph{Proof of Corollary \ref{KMTforcocycles}.} 
It follows from Theorem 2 (ii) in  \cite{CDJ} that  $n^{-1} \E_{\nu} \big ( (S_n - n \lambda_{\mu} )^2\big ) \rightarrow \sigma^2$ as $n \rightarrow \infty$. Then, the strong invariance principle  in stationary regime, i.e. for $(S_n)_{n\ge 1}$,  is a direct application of Theorem  \ref{thmsubexpo} by taking into account the estimate given in Proposition \ref{propsubexpodecreasecoeff} below. 
The proof of Proposition \ref{propsubexpodecreasecoeff} is given in appendix.

\begin{Proposition}  \label{propsubexpodecreasecoeff}
Assume that $\mu$ satisfies \eqref{condsubexpo} with $\gamma >0$.  Then, there exists $\beta >0$ such that, for any positive integer $k$, 
\begin{equation}\label{coef-expo}
 \sup_{ {\bar x}, {\bar y}  \in X }  \BBE (  |X_{k, {\bar x}}- X_{k, {\bar y}}|   )  \ll   {\rm e}^{- \beta  k^{ {\rm min} ( \gamma ,1)}}  \, .
\end{equation}
\end{Proposition}

\medskip

 As in \cite{CDJ} we obtain the result for the matrix norms $(\log(\|A_n\|)_{n\ge 1}$ from the stationary regime, 
using an argument from \cite{BL}. Indeed, it suffices to apply the 
estimate given just after equation (66) of \cite{CDJ} to get the 
ASIP in that case. Similarly, we get the ASIP for $(S_{n,\bar x})_{n\ge 1}$, for any $x\in S^{d-1}$. 

\medskip

It remains to handle the case of the matrix coefficients and of the spectral radius. The proof follows a well-known scheme that has been used in \cite{BQ-book} under exponential moments and in 
\cite{BQ} and \cite{CDMP} under polynomial moments. Those proofs rely on large deviations inequalities. In the case of sub-exponential moments, the needed large deviations estimates were obtained in \cite{CDM-deviation}.

\medskip

 We start with an auxiliary result of independent interest, see Theorem 2.6 of \cite{GLX}, for a related result. 

\begin{Proposition}\label{regularity}
Under the condition of Corollary \ref{KMTforcocycles}, there exists 
$\eta>0$ such that 
$$
\sup_{\bar x\in X}\int_X {\rm e}^{\eta \left |\log \delta(\bar x, \bar y\rangle \right |^\gamma}d\nu(\bar y)<\infty\, ,
$$
where $\delta (\bar x, \bar y)=|\langle x,y\rangle| /(\|x\|\|y\|)$.
\end{Proposition}

We can now prove the result for the matrix coefficients. 

\medskip
Let $\eta$ be as in Proposition \ref{regularity}. Let $n\in \N$ and $\bar x\in X$. Using that $\nu$ is $\mu$-invariant and Markov's inequality, we have
\begin{align*}
\P\left(-\log\delta (W_n, \bar x)\ge (2\log n/\eta)^{1/\gamma}\right ) 
=\nu\left(\left \{\bar y\in X\, :\, {\rm e}^{\eta \left | \log \delta (\bar y, \bar x) \right |^\gamma}\ge n^2 \right \} \right )\le K/n^2\, .
\end{align*} 
Hence, by the Borel-Cantelli lemma we see that, $\P$-a.s. 
$$
-\log \delta (W_n, \bar x) \le (2\log n/\eta)^{1/\gamma}\, ,
$$
for all but finitely many $n\in \N$. Combining this with the computations after (67) (and page 1864) in \cite{CDJ}, we infer that 
for every $\bar x, \bar y\in X$, 
$$
\left |\log \delta (A_n \cdot \bar x, \bar y)\right |= O ((\log n)^{1/\gamma})\, ,
$$
and the desired result follows since 
$$
\log |\langle A_n x, y\rangle |= \log \|A_n x\|+\log \|y\| 
+\log \delta (A_n \cdot \bar x, \bar y) \, .$$

\medskip

It remains to prove the ASIP for the spectral radius.  The result follows from  the following lemma, with 
$\ell \ge (2\log n /c)^{1/\gamma}$, combined with the Borel-Cantelli lemma.

\begin{Lemma}
Under the assumption of Corollary \ref{KMTforcocycles}, for every $\varepsilon>0$, there exists $C, c>0$ and $\ell_0\in \N$ such that 
for every $\ell_0\le \ell \le n$, 
$$
\P(\log \rho (A_n)-\|A_n\|\ge -\varepsilon \ell)\ge 1- C {\rm e}^{-c\ell^\gamma}\, .
$$ 
\end{Lemma}

The lemma  may be proved as 
Lemma 8 in \cite{CDMP} using our Lemma \ref{large-deviations} from the Appendix instead of Lemma 6 in \cite{CDMP}. 

\subsection{Products of iid positive matrices}

Let us assume now  that $G$  is the semi-group of 
$d$-dimensional positive allowable matrices:  by positive, we mean that all entries are greater than or equal to 0, by allowable, we mean that  any row and any column admits a strictly positive element. 

\medskip
In this case, $G$ is acting on  $X=S^{d-1}\cap (\R^+)^d$ as follows: 
$g\cdot x=gx/\|gx\|$, and we define $N(g)$ by $N(g)=\max(\|g\|, 1/v(g))$, where $v(g):= \inf_{x\in X}\|gx\|$. Here, $\|\cdot \|$ stands for the $\ell^1$ norm on $\R^d$ which is more convenient for the problem, see \cite{CDM23}. This change does not affect the final results since all norms on $\R^d$ are equivalent.

\medskip

Then,  similarly to the previous section, we define polynomial moments of order $p\ge 1$ for $\mu$ as well as sub-exponential 
moments of order $\gamma\in ]0,1]$.

\medskip

We shall also say that $\mu$ is strictly contracting if its support contains a matrix whose all entries are \emph{strictly} positive. 
When $\mu$ is stricty contracting there exists a unique $\mu$-invariant probability $\nu$ on ${\mathcal B}(X)$, see Section 3 of \cite{CDM23}.

\medskip

We have the following analogue to Corollary \ref{KMTforcocycles} with $S_{n,x}= \log \|A_nx \|$, $x\in X$, where $A_n = \varepsilon_n \cdots \varepsilon_1$ with $(\varepsilon_n)_{n \geq 1}$  a sequence of  independent random matrices taking values in $G$, $d \geq 2$, with common distribution $\mu$. 
Let also $W_0$ be an
$X$-valued  random variable with law $\nu$ and 
 independent of $(\varepsilon_n)_{n \geq 1}$, and $S_n=\log \|A_nW_0\|$. The following corollary was announced in \cite{CDM23}.

 \begin{corollary} \label{KMT-positive} Let $\mu$ be a strictly contracting  probability measure on ${\mathcal B}(G)$. Assume that $\mu$ has a sub-exponential moment of order $\gamma \in ]0,1]$. Then $n^{-1} \E_{\nu} \big ( (S_n - n \lambda_{\mu} )^2\big ) \rightarrow \sigma^2$ as $n \rightarrow \infty$ and for every (fixed) $x\in X$, one can redefine $(S_{n,x})_{n \geq 1}$ without changing its distribution on a (richer) probability space on which 
there exist iid random variables $(N_i)_{i \geq 1}$ with common distribution ${\mathcal N} (0, \sigma^2)$, such that,
\begin{equation*} \label{asip-RW}
S_{n,{ x}} - n \lambda_{\mu} - \sum_{i=1}^n N_i = O( (\log n)^{2+1/\gamma}) \, \text{a.s.}
\end{equation*}
Moreover,the result remains true  if we replace the sequence $(S_{n, x})_{n\ge 1}$ 
with $(S_n)_{n\ge 1}$, $(\log \|A_n\|)_{n\ge 1}$, 
$(\log \rho(A_n))_{n\ge 1}$, $(\log v(A_n))_{n\ge 1}$ or $(\log |\langle A_n x, 
y\rangle |)_{n\ge 1}$,  for some $x,y\in X$.
\end{corollary}
\begin{remark} When $\mu$ is aperiodic, see Definition 5.1 of \cite{CDM23}, then, by Proposition 5.2 of \cite{CDM23}, $\sigma^2>0$. When $\mu$ admits super-exponential moments, it is 
possible to improve the rate as in the remark after Corollary \ref{KMTforcocycles}.
\end{remark}
\medskip

\noindent \emph{Proof of Corollary \ref{KMT-positive}.} Again, we start with the stationary regime. The existence of the asymptotic variance follows from Proposition 5.2 of \cite{CDM23}.

\medskip

By Proposition 3.2 of \cite{CDM23}, when  $\mu$  is strictly contracting and admits a moment of order  greater than one, \eqref{coef-expo} holds. Then, the result (for $(S_n)_{n\ge 1} $) follows from Theorem \ref{thmsubexpo}.

\medskip

 The ASIP for $(\log \|A_n\|)_{n\ge 1}$ and  $(\log v(A_n))_{n\ge 1}$ 
 follows from (4.3) in \cite{CDM23} and the fact that for every 
 $x\in X$, $v(A_n)\le \|A_nx\|\le \|A_n\|$. Similary, we obtain  the ASIP for $(S_{n,x})_{n\ge 1}$ for any $x\in X$. 
 
 \medskip
 
 Then, the ASIP for $(\log \rho(A_n))_{n\ge 1}$ follows from (4.3) in 
 \cite{CDM23}, combined with the inequalities, for every $n\ge 1$,
 $\log v(A_n)\le \log \rho(A_n)\le \log \|A_n\|$.

 \medskip
 
 The ASIP for the matrix coefficients is proved in \cite{CDM23}, see Theorem 6.4. 

\medskip

\subsection{Lipschitz autoregressive models} 

We consider the autoregressive Lipschitz model introduced by  Dedecker and Rio \cite{DR00}.
Let $\tau\in[0,1)$, $C\in (0,1]$ and $f\, :\, \R\to \R$ a $1$-Lipschitz 
function such that 

$$
f(0)=0\quad \mbox{ and } \quad |f'(t)|\le 1-\frac{C}{(1+|t|)^\tau} 
\quad \mbox{ for almost every $t$} \, .
$$
Let $(\varepsilon_i)_{i\ge 1}$ be iid real-valued random valued with common law $\mu$ and define for any $n \geq 1$
\beq \label{defLAM}
W_n = f (W_{n-1} ) + \varepsilon_n \, , \, \mbox{ with $W_0$ independent of $(\varepsilon_i)_{i\ge 1}$.}
\eeq
Let $S_n (g) = \sum_{k=1}^n g(W_i)$ for any measurable function $g$.  Assume that $\mu$ has a sub-exponential moment of order $\eta \in ]0,1]$ (hence 
$\BBE ( {\rm e}^{c |\varepsilon_0|^{\eta} } ) < \infty$, for some $c >0$). For this model, there exists an unique invariant probability measure $\nu$ 
(see Proposition 2 of \cite{DR00}) and the following result holds: 

%
%
%
%
%


\begin{corollary}
Assume that $\mu$ has a sub-exponential moment of order $\eta \in ]0,1]$. Let $(W_n)_{n \geq 0}$ be defined by \eqref{defLAM} with ${\mathcal L} (W_0) = \nu$. 
Assume furthermore that $g$ is Lipschitz and that, for some $\zeta \in [0,1]$ and some positive constant $\kappa$, $|g(x) | \leq \kappa ( 1 + |x|^\zeta)$. Suppose that $\tau+ \zeta >0$. Then,  $n^{-1}{\rm Var} (S_n(g) ) \rightarrow \sigma^2 (g)$ as $n \rightarrow \infty$ and   
one can redefine $(W_n)_{n \geq 0}$ without changing its distribution on a (richer) probability space on which 
there exist iid random variables $(N_i)_{i \geq 1}$ with common distribution ${\mathcal N} (0, \sigma^2(g))$, such that,
\[
S_n(g)-n \nu(g)  -  \sum_{i=1}^n N_i = O((\log n)^{1+1/\gamma_1+1/\gamma_2}) \, \text{ ${\mathbb P}$-a.s.}
\]
where $\gamma_1 = \eta (1- \tau) ( \eta (1- \tau)  + \tau )^{-1}$ and $\gamma_2 = \eta (1- \tau)  / \zeta$. 
\end{corollary}
\noindent {\bf Proof.}  The result comes from an application of Theorem \ref{thmsubexpo} by taking into account Proposition 3 of \cite{MPR} and its proof. More precisely,  
\cite[(4.3)]{MPR} and \cite[Remark 5]{MPR}  ensure that $g(W_1)-\nu(g)$ satisfies \eqref{condh} with index $\gamma_2 = \eta (1- \tau)  / \zeta$. In addition, as it appears clearly in the proof of \cite[(4.4)]{MPR},  the control of the $\tau$-coefficients is done via a control of the coupling coefficients $\delta (n)$ (see \cite[(4.6)]{MPR}).    $\square$

\section{Appendix} 

\subsection{Proof of Proposition \ref{propsubexpodecreasecoeff}} 

The proof of the proposition is based on the following lemma whose statement needs the following notation: For  ${\bar x}, {\bar y}  \in X$, define\[
d( {\bar x}, {\bar y} ) := \frac{ \Vert  x \wedge y\Vert }{\Vert x \Vert \Vert y \Vert }
\]
where $\wedge$  stands for the exterior product, see e.g.  \cite[page 61]{BL}  for the definition and some
properties.  Then $d$ is a metric on $X$.
\begin{Lemma}  \label{lmasubexpodecrease}
Assume that $\mu$ satisfies \eqref{condsubexpo} with $\gamma >0$. Then,  for any positive integer $k$ and any $\kappa >0$,  there exists $\ell >0$, not depending on $k$, such that 
\[
\max_{k \leq j \leq 2k} \sup_{ {\bar x}, {\bar y}  \in X , {\bar x} \neq  {\bar y}} {\mathbb P}  \big (  \log ( d (A_{j-1} \cdot {\bar x} , A_{j-1} \cdot {\bar y} ))  \geq - \ell k  \big )   \ll   {\rm e}^{- \kappa k^{ {\rm min} ( \gamma ,1)}} \, .
\]
\end{Lemma}
{\it Proof of Lemma \ref{lmasubexpodecrease}.} Without loss of generality, let us assume that $c >4$ in condition \eqref{condsubexpo}. The general case can be deduced by normalizing suitably the quantities of interest. We proceed with similar arguments as in the proof of Lemma 6 in \cite{CDJ}; namely, we use a martingale decomposition. 

Since $\int_G \log N(g)  \mu (d g)  < \infty$,  we  can define the following bounded function $F_1$:
\[
F_1 ( {\bar x}, {\bar y}  ) =  \int_G  \log ( d (g \cdot {\bar x} ,g \cdot {\bar y} ) / d  ( {\bar x}, {\bar y}  ) )  \mu (d g)   \, , \, \forall {\bar x}, {\bar y}  \in X , {\bar x} \neq  {\bar y} \, .
\]
Then we define the following centered cocycle:
\[
\sigma_1 ( g,  ( {\bar x}, {\bar y}  )) =  \log ( d (g \cdot {\bar x} , d (g \cdot {\bar y} ) / d  ( {\bar x}, {\bar y}  ) )  - F_1 ( {\bar x}, {\bar y}  ) \, .
\]
Setting 
\[
M_n =M_n  ( {\bar x}, {\bar y}  ) := \sum_{k=1}^n \sigma_1 ( \varepsilon_k,  ( A_{k-1} \cdot {\bar x}, A_{k-1} \cdot {\bar y}  ))  \mbox{ and } R_n = R_n   ( {\bar x}, {\bar y}  ) : = \sum_{k=1}^n F_1 (A_{k-1} \cdot {\bar x}, A_{k-1} \cdot {\bar y}  ) \, , 
\]
we have
\[
 \log ( d ( A_{n} \cdot {\bar x}, A_{n} \cdot {\bar y}  ) / d  ( {\bar x}, {\bar y}  ) ) = M_n +R_n \, .
\]
Hence, since $d  ( {\bar x}, {\bar y}  ) \leq 1$,  the lemma will be proved if one can show that for any $\kappa >0$, there exists a $\ell >0$ such that 
\beq \label{lmasubReste}
\max_{k \leq j \leq 2k} \sup_{ {\bar x}, {\bar y}  \in X , {\bar x} \neq  {\bar y}} {\mathbb P}  \big ( R_j   ( {\bar x}, {\bar y}  )  \geq - 2 \ell k  \big )  \ll    {\rm e}^{- \kappa k^{{\rm min} ( \gamma,1)}}
\eeq
and
\beq \label{lmasubMartingale}
\max_{k \leq j \leq 2k} \sup_{ {\bar x}, {\bar y}  \in X , {\bar x} \neq  {\bar y}} {\mathbb P}  \big (  | M_j   ( {\bar x}, {\bar y}  ) |   \geq  \ell k  \big ) \ll    {\rm e}^{- \kappa k^{{\rm min} ( \gamma,1)}} \, .
\eeq
The estimate  \eqref{lmasubReste} comes from the proof of (53) in \cite{CDJ}. Indeed from the last inequality in the proof of \cite[(53), page 1857]{CDJ}, there exist $\rho \in ]0,1[$, $n_0 \in {\mathbb N}$ and $\alpha >0$ such that 
\[
\max_{k \leq j \leq 2k} \sup_{ {\bar x}, {\bar y}  \in X , {\bar x} \neq  {\bar y}} {\mathbb P}  \big ( R_j   ( {\bar x}, {\bar y}  )  \geq - \alpha k  \big )  
\leq C \rho^{k/(2n_0)} \, .
\]
Hence \eqref{lmasubReste} follows by taking $\ell \geq \alpha n_0 / (|\log \rho|)$. 

 We turn to the proof of \eqref{lmasubMartingale}.  Let 
$D_k= \sigma_1 ( \varepsilon_k,  ( A_{k-1} \cdot {\bar x}, A_{k-1} \cdot {\bar y}  )) $ and ${\mathcal F}_k= \sigma (W_0, \varepsilon_1, \ldots, \varepsilon_k )$. 
Note that  $(D_k, {\mathcal F}_k)_{k>0}$ is a sequence of martingale differences. Moreover setting  $X = 4 \log (N(g))$,  Lemma 4.2 page 56  in \cite{BL} (see also their Proposition 5.3 page 62) ensures that 
\[
\Vert {\mathbb P} ( D_k   \geq x   | {\mathcal F_{k-1}}  )  \Vert_{\infty} \leq  {\mathbb P} ( X + \BBE (X)    \geq x )  \, ,
\]
implying that 
\[
 \Big \Vert \sum_{k=1}^n \BBE ( D^2_k  \exp \{  | D_k |^{\gamma} \} | {\mathcal F_{k-1}} )  \Big \Vert_{\infty} \leq n  \BBE ( (X + \BBE (X)  ) ^2 \exp \{ (X+ \BBE(X) )^{\gamma} \} )  :=n  K \, .
\] 
Note that since $\mu$ satisfies \eqref{condsubexpo} with $c >4$, we have $ K< \infty$. It follows from \cite[Theorem 2.1]{FGL} (in case $\gamma \in ]0,1[$) or from 
\cite[Theorem 1.1]{LW} (when $\gamma \geq 1$) that there exist $c_1$ and $c_2$ such that, for any $k \geq 1$ and any $\ell >0$, 
\[
\max_{k \leq j \leq 2k} \sup_{ {\bar x}, {\bar y}  \in X , {\bar x} \neq  {\bar y}} {\mathbb P}  \big (  | M_j   ( {\bar x}, {\bar y}  ) |   \geq  \ell k  \big ) 
\ll \exp ( - c_1 \ell^2 k/ K ) +  \exp ( - c_2 ( \ell k)^{{\rm min}(\gamma,1)} ) \, ,
\]
which proves \eqref{lmasubMartingale} and ends the proof of the lemma. $\square$

\smallskip

We turn now to the proof of Proposition \ref{propsubexpodecreasecoeff}. 

\smallskip

\noindent {\it Proof of Proposition \ref{propsubexpodecreasecoeff}.}
Let ${\bar x}, {\bar y}  \in X $ and $\kappa >0$. Let $\ell >0$ be as in Lemma  \ref{lmasubexpodecrease}.  We start from inequality (18) in \cite{CDJ}, namely: setting 
$$
A= \{ \log N( \varepsilon_k ) \geq  \ell k /2 \} \quad  \text{and} \quad B= \{   \log ( d (A_{k-1} \cdot {\bar x} , d (A_{k-1} \cdot {\bar y} ))  \geq - \ell k   \},$$ we write 
\begin{multline*}
|X_{k, {\bar x}}- X_{k, {\bar y}}| \leq | \sigma ( \varepsilon_k , A_{k-1} \cdot {\bar x} ) -   \sigma ( \varepsilon_k , A_{k-1} \cdot {\bar y} )  | {\bf 1}_A +  | \sigma ( \varepsilon_k , A_{k-1} \cdot {\bar x} ) -   \sigma ( \varepsilon_k , A_{k-1} \cdot {\bar y} )  | {\bf 1}_B \\
+  | \sigma ( \varepsilon_k , A_{k-1} \cdot {\bar x} ) -   \sigma ( \varepsilon_k , A_{k-1} \cdot {\bar y} )  | {\bf 1}_{\{A^c \cap B^c \} } \, .
\end{multline*}
Recall the estimates (13) and (14) in \cite{CDJ}: there exists $C>0$ such that for every ${\bar x}, {\bar y}  \in X $,
\[
 |  \sigma (g , {\bar x} )  -  \sigma (g ,  {\bar y} )  | \leq C N(g) d (  {\bar x},  {\bar y} ) \, , 
\]
and
\[
 |  \sigma (g , {\bar x} )  | \leq \log N(g)  \, .
\]
We then infer that 
\begin{multline} \label{ECP1}
\BBE ( |X_{k, {\bar x}}- X_{k, {\bar y}}|  ) \leq  2 \BBE (  ( \log N( \varepsilon_k ) )  {\bf 1}_A )  + 2 \BBE (  ( \log N( \varepsilon_k )  )  {\bf 1}_B ) \\
+  C  \BBE  \big ( N( \varepsilon_k )   d (  A_{k-1} \cdot {\bar x}  ,   A_{k-1} \cdot {\bar y}  ) | {\bf 1}_{\{A^c \cap B^c \} } \big )   \, .
\end{multline}
Now
\beq  \label{ECP2}
 \BBE (  (\log N( \varepsilon_k ) )  {\bf 1}_A )  \leq {\rm e}^{- (  \ell  k /2 )^{\gamma}} \BBE \big (   \log N(g)   {\rm e}^{ (\log N (g))^{\gamma}}   \big )  \, .
\eeq
On another hand, by independence between $\varepsilon_k$ and $B$, and Lemma \ref{lmasubexpodecrease}, 
\beq  \label{ECP3}
 \BBE (  (\log N( \varepsilon_k ) )  {\bf 1}_B )  = {\mathbb P} (B) \BBE (  \log N( g) )    \leq {\rm e}^{- \kappa k^{\gamma}} \BBE (  \log N( g) )   \, .
\eeq
Moreover,
\beq  \label{ECP4}
  \BBE  \big ( N( \varepsilon_k )   d (  A_{k-1} \cdot {\bar x}  ,   A_{k-1} \cdot {\bar y}  ) | {\bf 1}_{\{A^c \cap B^c \} } \big )  \leq 
  {\rm e}^{- \ell k }   \BBE  \big ( N(  g) {\bf 1}_{\{   \log N( g ) \leq \ell k /2 \}}   \big )  \leq    {\rm e}^{- \ell k /2 }  
 \, .
\eeq
Starting from \eqref{ECP1} and taking into account the upper bounds  \eqref{ECP2}-\eqref{ECP4}, the lemma follows. $\square$

\subsection{Proof of Proposition \ref{regularity}}
 It suffices to find $\eta>0$ such that (notice that $\delta ( \cdot, \cdot )\le 1$)
$$
\sum_{n\ge 1} {\rm e}^{ n^\gamma}\sup_{\bar x\in X}\nu(\{\bar y\in X\, :\, - \log \delta(\bar x, \bar y ) \ge \eta n\})<\infty\, .
$$
Now, using that $\nu$ is $\mu$-invariant, for every $n\in \N$ and $\bar x\in X$, 
\begin{align*}
\nu(\{\bar y\in X\, :\, - \log \delta(\bar x, \bar y ) \ge \eta n \})
  &  =\nu\otimes \mu^{*n}(\{(\bar y,g)\in X\times G\, :\, - \log \delta(\bar x, g\cdot\bar y )\ge \eta n\})\\
    &  \le \sup_{\bar x,\bar y\in X} \P ( \log \delta(\bar x, A_n \cdot\bar y ) \le -\eta n) \, .
\end{align*}
The fact that one may find such a $\eta>0$ may be proved as Lemma 
4.8 (and more particularly (4.13)) of \cite{BQ}, using the 
next lemma. 

\begin{Lemma}\label{large-deviations}
Under the assumption of Corollary \ref{KMTforcocycles}, for every 
$\varepsilon>0$, there exist $C, c>0$, such that 
\begin{gather}
\label{BQ1}\sup_{x\in S^{d-1}} \P \left(\max_{1 \leq k \leq n}|\log \|A_k x\|-k \lambda_\mu|>\varepsilon n \right)
\le C{\rm e}^{-c n^\gamma}\, ,\\
\label{BQ2} \P\left(\max_{1 \leq k \leq n}|\log \|A_k\|-n \lambda_\mu|>\varepsilon n \right)
\le C{\rm e}^{-c n^\gamma}\, ;\\
\label{BQ3} \P\left (\max_{1 \leq k \leq n}|\log \|\Lambda^2(A_k)\|-k (\lambda_\mu 
 +\gamma_\mu)|>\varepsilon n \right)
\le C{\rm e}^{-c n^\gamma}\, .
\end{gather}
\end{Lemma}
\begin{remark} Let us recall that, for any $A \in GL_d({\mathbb R})$, $\Lambda^2(A)$ is the matrix on $\Lambda^2({\mathbb R }^d)$ defined by $\Lambda^2(A)(x \wedge y) = Ax \wedge A y$, where $\wedge$ is the exterior product.  In addition, in  \eqref{BQ3} $,\gamma_\mu$ is the second Lyapunov exponent of $\mu$. With the notations of 
\cite[Section 14]{BQ-book}, $\lambda_\mu$ is denoted either $\lambda_{1,\mu}$ or $\lambda_1$, while $\gamma_\mu$ is either denoted $\lambda_{2,\mu}$ or $\lambda_2$.  
\end{remark}

\noindent \emph{Proof of Lemma \ref{large-deviations}}. \eqref{BQ1} follows directly from Theorem 2.2 of \cite{CDM-deviation}. Then 
\eqref{BQ2} follows from \eqref{BQ1} and  the fact that there exists 
$C>0$ such that for every $g\in G$ and every $x\in S^{d-1}$, 
$$
\|gx\|\le \|g\|\le C\max_{1\le i\le d}\|g e_i\|\, ,
$$
with $\{e_i\, :\, 1\le i\le d\}$ the canonical basis of $\R^d$.

\medskip

The proof of \eqref{BQ3} may be done similarly as above. The group $G$ acts on $P(\wedge^2(\R^d))$, the projective space of 
$\wedge^2 (\R^d)$, by $g\cdot \bar \xi:=\wedge^2 (g) \xi  /\|\xi\|$, for every $g\in G$ and $\xi \in \wedge^2 (\R^d)$, where $\bar \xi$ is the projection of $\xi$ on $P(\wedge^2 (\R^d))$. Then, 
$(g,\bar \xi)\mapsto \|\wedge^2(g)(\xi) \|/\|\xi\|$ defines a cocycle and one may 
prove, similarly to \eqref{BQ1} that (see Section 6 of \cite{CDM-deviation} where the situation of general cocycles is considered)
$$
\sup_{\xi\in \wedge ^2(\R^d)\, :\, \|\xi\|=1}\P\left (\max_{1\le k\le n} |
\log \|\wedge^2(A_k)(\xi)-k(\lambda_\mu+\gamma_\mu)|> \varepsilon n \right )
\le C{\rm e}^{-cn^\gamma}\, .
$$
Then, we deduce \eqref{BQ3} as we deduced \eqref{BQ2} from \eqref{BQ1}.


\medskip




\begin{thebibliography}{99}

%

\bibitem{BQ} Benoist, Y. and  Quint J.-F. \emph{Central limit theorem for linear groups}, Ann. Probab.  44  (2016),  no. 2, 1308-1340.

\bibitem{BQ-book}  Benoist, Y. and  Quint J.-F. Random walks on reductive groups, Ergebnisse der Mathematik und ihrer Grenzgebiete. 3. Folge., vol. 62, Springer, 2016.

\bibitem{BLW14} Berkes, I.,  Liu, W. and  Wu, W. B. 
Koml\'os-Major-Tusn\'ady approximation under dependence. {\it Ann. Probab.} {\bf 42} (2014), no. 2, 794-817.

\bibitem{BL} Bougerol, P. and Lacroix, J. Products of random matrices with applications to Schr\"odinger operators. 
Progress in Probability and Statistics, 8. Birkh\"auser Boston, Inc., Boston, MA, 1985

%
%

\bibitem{CDJ} Cuny, C.,  Dedecker, J. and Jan, C. Limit theorems for the left random walk on $GL_d({\mathbb R})$. \textit{Ann. Inst. Henri Poincaré Probab. Stat.} {\bf 53} (2017), no. 4, 1839--1865. 

\bibitem{CDKM0} Cuny, C., Dedecker, J., Korepanov, A. and  Merlev\`ede, F. Rates in almost sure invariance principle for slowly mixing dynamical systems. {\it Ergodic Theory Dynam. Systems} {\bf 40} (2020), no. 9, 2317--2348.


\bibitem{CDKM}  Cuny, C., Dedecker, J., Korepanov, A. and  Merlev\`ede, F. Rates in almost sure invariance principle for quickly mixing dynamical systems. {\it Stoch. Dyn.} {\bf 20} (2020), no. 1, 2050002, 28 pp. 

\bibitem{CDM-deviation} Cuny, C., Dedecker, J. and Merlev\`ede, F., Large and moderate deviations for the left random walk on $GL_d({\mathbb R})$. ALEA Lat. Am. J. Probab. Math. Stat. 14 (2017), no. 1, 503–527.

\bibitem{CDM} Cuny, C.,  Dedecker, J. and Merlev\`ede, F. On the Koml\'os, Major and Tusn\'ady strong approximation for some classes of random iterates. {\it Stochastic Process. Appl.} {\bf 128} (2018), no. 4, 1347--1385. 

\bibitem{CDMbis} Cuny, C., Dedecker, J. and  Merlevède, F. An alternative to the coupling of Berkes-Liu-Wu for strong approximations. {\it Chaos Solitons Fractals} {\bf 106} (2018), 233--242. 

\bibitem{CDM23} Cuny, C.,  Dedecker, J. and Merlev\`ede, F. Limit theorems for iid products of positive matrices. \texttt{hal-03923713}

\bibitem{CDMP} Cuny, C., Dedecker, J., Merlev\`ede, F. and  Peligrad, M., Berry-Esseen type bounds for the matrix coefficients and the spectral radius of the left random walk on $GL_d({\mathbb R})$. C. R. Math. Acad. Sci. Paris 360 (2022), 475-482.

%
%

\bibitem{DD03}  Dedecker, J. and  Doukhan, P.  A new covariance inequality and applications. {\it Stochastic Process. Appl.} {\bf 106} (2003), no. 1, 63-80.

%


\bibitem{DR00}  Dedecker, J. and  Rio, E.  On the functional central limit theorem for stationary processes. {\it Ann. Inst. H. Poincaré Probab. Statist.} {\bf 36} (2000), no. 1, 1-34.

 \bibitem{FGL} Fan, X.,  Grama, I. and  Liu, Q.  Deviation inequalities for martingales with applications. {\it J. Math. Anal. Appl.}  {\bf 448} (2017), no. 1, 538--566.

\bibitem{FK} Furstenberg, H. and  Kesten, H. Products of Random Matrices. {\it Ann. Math. Statist.} {\bf 31} (1960), no. 2, 457-469.


\bibitem{GLX} Grama, I., Liu Q. and Xiao H., \emph{Limit theorems on the coefficients of random walks on the general linear group}, 
\texttt{arXiv:2111.10569}

%
%

\bibitem{Kor} Korepanov, Alexey Rates in almost sure invariance principle for dynamical systems with some hyperbolicity. {\it Comm. Math. Phys.} {\bf 363} (2018), no. 1, 173--190. 

\bibitem{LW} Liu, Q. and  Watbled, F.  Exponential inequalities for martingales and asymptotic properties of the free energy of directed polymers in a random environment. {\it Stochastic Process. Appl.} {\bf 119} (2009), no. 10, 3101--3132.  

%
%
%

\bibitem{MPR} Merlev\`ede, F.,  Peligrad, M. and  Rio, E.  A Bernstein type inequality and moderate deviations for weakly dependent sequences. \textit{Probab. Theory Related Fields} {\bf 151} (2011), no. 3-4, 435-474. 
 
 
%
%
%


\bibitem{Ri00} Rio, E.  {\it Th\'eorie asymptotique des processus al\'eatoires faiblement dépendants}, Math\'ematiques $\&$ Applications (Berlin), vol. 31, Springer-Verlag, Berlin,
2000.

%
%


\end{thebibliography}
\end{document}